\RequirePackage{fix-cm}
\documentclass[smallextended]{svjour3}       
\smartqed  
\usepackage{graphicx}

\begin{document}

\title{ Computing the $2-$adic complexity of two classes of  Ding-Helleseth generalized cyclotomic sequences of period of twin prime products\thanks{The work is financially supported by National Natural Science Foundation of China (No. 61902429, No.11775306), Shandong Provincial Natural Science Foundation of China (No. ZR2017MA001, ZR2019MF070), Fundamental Research Funds for the Central Universities (No. 19CX02058A, No.17CX02030A), the Open Research Fund from Shandong provincial Key Laboratory of Computer
Networks, Grant No. SDKLCN-2018-02, Key Laboratory of Applied Mathematics of Fujian Province University (Putian University) (No. SX201702, No.SX201806), the Projects of International Cooperation and Exchanges NSFC-RFBR (No. 61911530130), and International Cooperation Exchange Fund of China University of Petroleum (UPCIEF2019020).}}



\author{Ming Yan  \and Tongjiang Yan  \and Yu Li}


\institute{
Ming Yan \at
College of Sciences,
China University of Petroleum,
Qingdao 266555,
Shandong, China\\
\email{yanming$\_a1$@163.com}
\and
Tongjiang Yan\at
College of Sciences,
China University of Petroleum,
Qingdao 266555,
Shandong, China\\
\email{yantoji@163.com}
(corresponding author)
\and
Yu Li \at
School of Cyber Security,
University of Chinese Academy of Sciences,
Beijing 100049,
China\\
\email{liyu1@iie.ac.cn}
}
\date{Received: date / Accepted: date}

\maketitle

\begin{abstract}
 This paper contributes to compute $2$-adic complexity of two classes of Ding-Helleseth generalized cyclotomic sequences. Results show that $2$-adic complexity of these sequences is good enough to resist the attack  by  the rational approximation algorithm.
\end{abstract}

\section{Introduction}
\label{section 1}
Pseudo-random sequences that can be widely used in cryptography  must have good security characteristics \cite{dingcunsheng1}. In 1995, $2$-adic complexity of sequences was introduced  as an important pseudo-randoman measure by A.Klapper, which is the length of the shortest feedback with carry shift register.  Its value cannot be less than half of the period, otherwise the sequence will be vulnerable to be attacked by $2$-adic of the rational approximation algorithm (RAA) \cite{klapper2}. A. Klapper, Tian and Qi showed that the $2$-adic complexity of all the binary $m$-sequences is maximal \cite{klapper3,tian-qi4}. Then Xiong $et$ $al$. and Hu proved that all the known sequences with ideal 2-level autocorrelation have maximum $2$-adic complexity by a new method using circulant matrices. Moreover, they also proved that the $2$-adic complexities of Legendre sequences, Ding-Helleseth-Lam sequences and two classes of interleaved sequences with optimal autocorrelation are also maximal \cite{huhonggang4-5,xionghai5,xionghai6}. Then, Sun and Yan $et$ $al$. give a lower bound on the $2$-adic complexity of the modified Jacobi sequence and the $2$-adic complexity of a class of binary sequences with almost optimal autocorrelation \cite{sunyuhua7,sunyuhua8}.

In 2017, Vladimir Edemskiy and Nikita Sokolovskiy have proved two classess of Ding-Helleseth generalized cyclotomic sequences of period $pq$ possess good linear complexity \cite{Vladimir9}. In this paper, we will give the $2$-adic complexity of these sequences.
\section{Preliminaries}
Let $p$ and $q$ ($p<q$) be two odd distinct primes with $\gcd(p-1, q-1)=2$, $\frac{(p-1)(q-1)}{2}=e$. Define $N=pq$.
By  Chinese Remainder Theorem, there exists a common primitive root $g$ of both $p$ and $q$, and an integer $x$ satisfying $x\equiv g\pmod{p}$, and $x\equiv1\pmod{q}$.
Ding-Helleseth generalized cyclotomic classes $D_i^{(N)}$ of order $2$ with respect to $p$ and $q$ are defined as
\begin{equation}\label{eqD0D1}
\begin{array}{ll}
D_0^{(N)}&=\{g^{2t}x^i: t=0, 1, \cdots,
\frac{e}{2}-1;\ i=0, 1\};\\
D_1^{(N)}&=gD_0^{(N)}.
\end{array}
\end{equation}
Clearly,
\[ D_1^{(N)}=gD_0^{(N)}, Z_N^{\ast}=\ D_0^{(N)}\cup D_1^{(N)}, \ D_0^{(N)}\cap D_1^{(N)}=\emptyset,\]
where $\emptyset$ denotes the empty set, $Z_N^{*}$ denotes the set of all invertible elements of the residue class ring $Z_N$.

Put $P=\{p, 2p, \cdots, (q-1)p\}$, $Q=\{q, 2q, \cdots, (p-1)q\}$, and $R=\{0\}$. Then $D_0^{(N)}\cup D_1^{(N)}\cup P\cup Q \cup R$ is a partition of $Z_N$. In the following we investigate the $2$-adic complexity of two types of sequences based on them.

The cyclotomic numbers to the generalized cyclotomic of order 2 are defined by
\begin{eqnarray}\label{eq(i,j)-n}
(i,j)_2^{(n)}&=&|(D_{i}^{(n)}+1)\cap D_{j}^{(n)}|,\ \mathrm{for\ all}\ i,j=0,1.
\end{eqnarray}

Assume $D_{0}^{(q)}$ and $D_{1}^{(q)}$ are the  quadratic residue class and  non-quadratic residue class of $Z_q$. Define the corresponding cyclotomic numbers
\begin{eqnarray}\label{eq(i,j)-q}
(i,j)_2^{(q)}=|(D_{i}^{(q)}+1)\cap D_{j}^{(q)}|,\ \mathrm{for\ all}\ i,j=0,1.
\end{eqnarray}

Let $\xi_N=e^{\frac{2\pi\sqrt{-1}}{N}}$ be a $N$th complex primitive root of unity. Then the additive character $\chi_N$ of $Z_N$ is given by
\begin{equation}\label{eqxi-n}
\chi_N(x)=\xi_N^{x},\ x\in Z_N,
\end{equation}
and Gauss periods of order $d$ are defined by
\begin{equation}\label{eqeta-n}
\eta_{i}^{(N)}=\sum_{x\in D_{i}^{(N)}}\chi_N(x),\ i=0,1.
\end{equation}
Correspondingly, Gauss periods for $Z_q$ are defined by
\begin{eqnarray}\label{eqeta-q}
\eta_{i}^{(q)}=\sum\limits_{x\in D_{i}^{(q)}}\chi_q(x),\ i=0,1,
\end{eqnarray}
where $\chi_q(x)=e^{\frac{2\pi\sqrt{-1}x}{q}}$ be the additive character of $Z_q$.

Let $N$ be a positive integer, $\{s_i\}_{i=0}^{N-1}$ be a binary sequence of period $N$, and
\begin{equation}\label{eqS(x)}
S(x)=\sum\limits_{i=0}^{N-1}s_{i}x^i\in Z [x].
\end{equation}
If we write
\begin{eqnarray}
\frac{S(2)}{2^N-1}=\frac{\sum\limits_{i=0}^{N-1}s_{i}2^i}{2^N-1}=\frac{m}{n},\ 0\leq m\leq n,\ \mathrm{gcd}(m,n)=1,\label{$2$-adic complexity}\nonumber
\end{eqnarray}
then the $2$-adic complexity $\Phi_{2}(s)$ of the sequence $\{s_i\}_{i=0}^{N-1}$ is defined as the integer $\lfloor\mathrm{log}_2n\rfloor$, i.e.,
\begin{equation}\label{eqphi-2}
\Phi_{2}(s)=\left\lfloor\mathrm{log}_2\frac{2^N-1}{\mathrm{gcd}(2^N-1,S(2))}\right\rfloor,\label{$2$-adic calculation}
\end{equation}
where $\lfloor x\rfloor$ is the greatest integer that is less than or equal to $x$.

Let $A=(a_{i,j})_{N\times N}$ be the circulant matrix defined by $a_{i,j}=s_{i-j\pmod{N}}$ and we view $A$ as a matrix  in the field  $C$ of complex numbers.

\section{Subsidiary lemmas}
For any $a\in Z_N$ and $B\subseteq Z_N$, we denote $aB=\{ab \mid b\in B\}$. Then we have the following properties.
\begin{lemma}\label{lmPQR}  \cite{sunyuhua7} Let the symbols be the same as before.
\begin{itemize}
\item[(1)]For each fixed $a\in D_i^{(N)}$, we have $aD_j^{(N)}=D_{(i+j)\pmod{2}}^{(N)}$, $aP=P$ and $aQ=Q$, where $i,j=0,1$.
\item[(2)] For each fixed $a\in P$, if $b$ runs through each element of $D_i^{(N)}$, $i=0,1$, then $ab$ runs  each element of $P$ exactly $\frac{p-1}{2}$ times. Symmetrically, for each fixed $a\in Q$, if $b$ runs through each element of $D_i^{(N)}$, $i=0,1$, then $ab$ runs each element of $Q$ exactly $\frac{q-1}{2}$ times.
\item[(3)]For each fixed $a\in P$, we have $aP=P$, $aQ=R$. Symmetrically, for each fixed $a\in Q$, we have $aQ=Q$, $aP=R$.
\end{itemize}
\end{lemma}
\begin{lemma}\label{lm-1}  \cite{wangqiuyan10}
Let the symbols be the same as before.
\begin{eqnarray}
-1\in\left\{
\begin{array}{ll}
D_0^{(N)},
&\ \textrm{if $q\equiv 1\pmod{4}$;}\nonumber\\
D_1^{(N)},
&\ \textrm{if $q\equiv 3\pmod{4}$.}\nonumber\\
\end{array}
\right.\label{two cases}
\end{eqnarray}
\end{lemma}
\begin{lemma}\label{lmcyc-num-q}  \cite{dingcunsheng1}
Let $q=2f+1$ be a prime, the cyclotic numbers of order 2 are given by

\noindent (1) $(0,0)_2^{(q)}=\frac{f-2}{2}; (0,1)_2^{(q)}=(1,0)_2^{(q)}=(1,1)_2^{(q)}=\frac{f}{2}$ if $f$ is even;

\noindent (2) $(0,0)_2^{(q)}=(1,0)_2^{(q)}=(1,1)_2^{(q)}=\frac{f-1}{2}; (0,1)_2^{(q)}=\frac{f+1}{2}$ otherwise.
\end{lemma}
\begin{lemma}\label{lmcyc-num-n}  \cite{Guass-sums11}
If $q\equiv1\pmod{4}$, then $$(1,1)_2^{(N)}=(1,0)_2^{(N)}=(0,1)_2^{(N)}=\frac{(p-2)(q-1)}{4}, (0,0)_2^{(N)}=\frac{(p-2)(q-5)}{4}.$$

If $q\equiv3\pmod{4}$, then $$(0,0)_2^{(N)}=(1,1)_2^{(N)}=(1,0)_2^{(N)}=\frac{(p-2)(q-3)}{4}, (0,1)_2^{(N)}=\frac{(p-2)(q+1)}{4}.$$
\end{lemma}
\begin{lemma}\label{lmy-x}  \cite{wangqiuyan10}
Let the symbols be the same as before.
\begin{eqnarray}
|(D_i^{(N)}+\omega)\cap D_j^{(N)}|=\left\{
\begin{array}{ll}
(p-1)(i,j)_2^{(q)},
&\ \textrm{if $\omega\in P$ and $(\frac{\omega}{q})=1$;}\nonumber\\
(p-1)(i+1,j+1)_2^{(q)},
&\ \textrm{if $\omega\in P$ and $(\frac{\omega}{q})=-1$;}\nonumber\\
\frac{(q-1)(p-2)}{2},
&\ \textrm{if $\omega\in Q$ and $i=j$;}\nonumber\\
0,
&\ \textrm{if $\omega\in Q$ and $i\neq j$,}\nonumber\\
\end{array}
\right.
\end{eqnarray}
where $(*,*)_2^{(q)}$ and $(\frac{*}{*})$ denote the cyclomic number of order 2 with respect to $q$ and the Legendre symbol, respectively.
\end{lemma}
\begin{lemma}\label{lmgcdA-N}\cite{sunyuhua7}
\noindent (1) $\mathrm{det}(A)=\prod_{t=0}^{N-1}S(\xi_N^t)$.

\noindent (2) If $\mathrm{det}(A)\neq0$, then
$\mathrm{gcd}\left(S(2),2^N-1\right) \mid \mathrm{gcd}\left(\mathrm{det}(A),2^N-1\right).$
\end{lemma}
\begin{lemma}\label{lmchi=1}  \cite{sunyuhua7}
It is well-known that
\begin{eqnarray}
\sum_{x\in P}\chi(x)=-1,\nonumber
\sum_{x\in Q}\chi(x)=-1,\nonumber
\sum_{i=0}^{1}\eta_i^{(N)}=1,\nonumber
\sum_{i=0}^{1}\eta_i^{(q)}=-1.\nonumber
\end{eqnarray}
\end{lemma}
\section{The $2$-adic complexity of the Ding-Helleseth generalized cyclotomic sequence of order 2}
\subsection{The $2$-adic complexity of the first class of sequences}
Define the Ding-Helleseth generalized cyclotomic sequence of order 2 with respect to primes $p$ and $q$ as
\begin{equation}\label{eqsi-p}
s_i=\left\{
\begin{array}{ll}
0, \ \mathrm{if}\ i\pmod{N}\in D_0^{(N)}\cup Q\cup R,\\
1, \ \mathrm{if}\ i\pmod{N}\in D_1^{(N)}\cup P.\\
\end{array}
\right.
\end{equation}
\begin{theorem}
Let $A$ be the circulant matrix based on the Ding-Helleseth generalized cyclotomic sequence defined by Equation (\ref{eqsi-p}), then
\begin{eqnarray}
&&\mathrm{det}(A)\nonumber\\
&=&\left\{
\begin{array}{llll}
(p+1)\left(\frac{q-1}{2}\right)^{p}\left(\left(p-1\right)^2\frac{q-1}{4}-p\right)^{\frac{q-1}{2}}\left(\frac{q-1}{4}\right)^{\frac{(p-1)(q-1)}{2}},
&\ \textrm{if $q\equiv 1\pmod{4}$;}\nonumber\\
(p+1)\left(\frac{q-1}{2}\right)^{p}\left(\left(p-1\right)^2\frac{q+1}{4}+p\right)^{\frac{q-1}{2}}\left(\frac{q+1}{4}\right)^{\frac{(p-1)(q-1)}{2}},
&\ \textrm{if $q\equiv 3\pmod{4}$.}\nonumber
\end{array}
\right.\label{two cases}
\end{eqnarray}
\end{theorem}
\noindent{\bf Proof.}
By Lemma \ref{lmgcdA-N},
\begin{eqnarray}
\mathrm{det}(A)&=&\prod_{a=0}^{N-1}S(\xi_N^a)\nonumber\\
&=&\prod_{a\in R}S(\xi_N^a)\prod_{a\in P}S(\xi_N^a)\prod_{a\in Q}S(\xi_N^a)\prod_{a\in D_0^{(N)}}S(\xi_N^a)\prod_{a\in D_1^{(N)}}S(\xi_N^a).\nonumber
\end{eqnarray}
By Equations (\ref{eqS(x)}) and (\ref{eqsi-p}),
\begin{eqnarray}
S(\xi_N^a)&=&\sum_{i\in D_1^{(N)}}\xi_N^{ai}+\sum_{i\in P}\xi_N^{ai}.\nonumber
\end{eqnarray}
(1) If $a\in R$, obviously,
\begin{eqnarray}
S(\xi_N^a)=|D_1^{(N)}|+|P|=\frac{(p+1)(q-1)}{2}.\nonumber
\end{eqnarray}
Thus$$\prod_{a\in R}S(\xi_N^a)=\frac{(p+1)(q-1)}{2}.$$
(2) If $a\in P$,
\begin{itemize}
\item[1)]  when $a\in D_0^{(q)}p$, by Lemmas \ref{lmPQR}, \ref{lmchi=1}, Equations (\ref{eqD0D1}) and (\ref{eqeta-q}),
\begin{eqnarray}
S(\xi_N^a)&=&(p-1)\eta_1^{(q)}-1;\nonumber
\end{eqnarray}
\item[2)]  when $a\in D_1^{(q)}p$, by Lemmas \ref{lmPQR}, \ref{lmchi=1}, Equations (\ref{eqD0D1}) and (\ref{eqeta-q}),
\begin{eqnarray}
S(\xi_N^a)&=&(p-1)\eta_0^{(q)}-1.\nonumber
\end{eqnarray}
\end{itemize}
And we have
\begin{eqnarray}
\prod_{a\in P}S(\xi_N^a)&=&\left(((p-1)\eta_0^{(q)}-1)((p-1)\eta_1^{(q)}-1)\right)^{\frac{q-1}{2}}\nonumber\\
&=&\left((p-1)^2\eta_0^{(q)}\eta_1^{(q)}-(p-1)(\eta_0^{(q)}+\eta_1^{(q)})+1\right)^{\frac{q-1}{2}}\nonumber\\
&=&\left((p-1)^2\eta_0^{(q)}\eta_1^{(q)}+p\right)^{\frac{q-1}{2}}.\nonumber
\end{eqnarray}
If $q\equiv 1\pmod{4}$, by Lemmas \ref{lm-1}, \ref{lmcyc-num-q}, \ref{lmchi=1} and Equation (\ref{eqeta-q}),
\begin{eqnarray}
\eta_0^{(q)}\eta_1^{(q)}&=&(0,1)_2^{(q)}\eta_0^{(q)}+(1,0)_2^{(q)}\eta_1^{(q)}\nonumber\\
&=&(0,1)_2^{(q)}(\eta_0^{(q)}+\eta_1^{(q)})\nonumber\\
&=&-\frac{q-1}{4}.\nonumber
\end{eqnarray}
If $q\equiv 3\pmod{4}$, by Lemmas \ref{lm-1}, \ref{lmcyc-num-q}, \ref{lmchi=1} and Equation (\ref{eqeta-q}),
\begin{eqnarray}
\eta_0^{(q)}\eta_1^{(q)}&=&|D_1^{(q)}|+(1,1)_2^{(q)}\eta_0^{(q)}+(0,0)_2^{(q)}\eta_1^{(q)}\nonumber\\
&=&\frac{q-1}{2}-\frac{q-3}{4}\nonumber\\
&=&\frac{q+1}{4}.\nonumber
\end{eqnarray}
So we get
\begin{eqnarray}
\prod_{a\in P}S(\xi_N^a)=\left\{
\begin{array}{ll}
\left((p-1)^2\frac{q-1}{4}-p\right)^{\frac{q-1}{2}},
&\ \textrm{if $q\equiv 1\pmod{4}$;}\nonumber\\
\left((p-1)^2\frac{q+1}{4}+p\right)^{\frac{q-1}{2}},
&\ \textrm{if $q\equiv 3\pmod{4}$.}\nonumber
\end{array}
\right.\label{two cases}
\end{eqnarray}
(3) If $a\in Q$, by Lemmas \ref{lmPQR} and \ref{lmchi=1},
\begin{eqnarray}
S(\xi_N^a)&=&\frac{q-1}{2}\sum_{k\in Q}\xi_N^k+(q-1)\nonumber\\
&=&\frac{q-1}{2}.\nonumber
\end{eqnarray}
Then$$\prod_{a\in Q}S(\xi_N^a)=\left(\frac{q-1}{2}\right)^{p-1}.$$
(4) If $a\in D_0^{(N)}$, by Lemmas \ref{lmPQR} and \ref{lmchi=1},
\begin{eqnarray}
S(\xi_N^a)&=&\eta_1^{(N)}-1=-\eta_0^{(N)}.\nonumber
\end{eqnarray}
Then$$\prod_{a\in D_0^{(N)}}S(\xi_N^a)=\left(\eta_0^{(N)}\right)^{\frac{(p-1)(q-1)}{2}}.$$
(5) Similarly, if $a\in D_1^{(N)}$, by Lemmas \ref{lmPQR} and \ref{lmchi=1},
\begin{eqnarray}
S(\xi_N^a)&=&\eta_0^{(N)}-1=-\eta_1^{(N)}.\nonumber
\end{eqnarray}
Then$$\prod_{a\in D_1^{(N)}}S(\xi_N^a)=\left(\eta_1^{(N)}\right)^{\frac{(p-1)(q-1)}{2}}.$$
Then we get
\begin{eqnarray}
&&\mathrm{det}(A)\nonumber\\
&=&\left\{
\begin{array}{ll}
&\frac{(p+1)(q-1)}{2}\left(\frac{q-1}{2}\right)^{p-1}\left(\left(p-1\right)^2\frac{q-1}{4}-p\right)^{\frac{q-1}{2}}\left
(\eta_0^{(N)}\eta_1^{(N)}\right)^{\frac{(p-1)(q-1)}{2}},\\
&\ \textrm{if $q\equiv 1 \pmod{4}$;}\nonumber\\
&\frac{(p+1)(q-1)}{2}\left(\frac{q-1}{2}\right)^{p-1}\left(\left(p-1\right)^2\frac{q+1}{4}+p\right)^{\frac{q-1}{2}}\left
(\eta_0^{(N)}\eta_1^{(N)}\right)^{\frac{(p-1)(q-1)}{2}},\\
&\ \textrm{if $q\equiv 3 \pmod{4}$.}\nonumber\\
\end{array}
\right.
\end{eqnarray}
When $q\equiv 1\pmod{4}$, $-1\in D_0^{(N)}$, by Lemmas \ref{lm-1} - \ref{lmy-x}, \ref{lmchi=1} and Equation (\ref{eqeta-n}), we have
\begin{eqnarray}
&&\eta_0^{(N)}\eta_1^{(N)}\nonumber\\
&=&\sum_{x\in D_0^{(N)}}\sum_{y\in D_1^{(N)}}\xi_N^{y-x}\nonumber\\
&=&(0,1)_2^{(N)}\eta_0^{(N)}+(1,0)_2^{(N)}\eta_1^{(N)}-(p-1)(0,1)_2^{(q)}\nonumber\\
&=&(0,1)_2^{(N)}-(p-1)(0,1)_2^{(q)}\nonumber\\
&=&-\frac{q-1}{4}.\nonumber
\end{eqnarray}
When $q\equiv 3\pmod{4}$, $-1\in D_1^{(N)}$, by Lemmas \ref{lm-1} - \ref{lmy-x}, \ref{lmchi=1} and Equation (\ref{eqeta-n}), we have
\begin{eqnarray}
&&\eta_0^{(N)}\eta_1^{(N)}\nonumber\\
&=&\sum_{x\in D_1^{(N)}}\sum_{y\in D_1^{(N)}}\xi_N^{y-x}\nonumber\\
&=&|D_1^{(N)}|+(1,1)_2^{(N)}\eta_0^{(N)}+(0,0)_2^{(N)}\eta_1^{(N)}-(p-1)(0,0)_2^{(q)}-\frac{(p-2)(q-1)}{2}\nonumber\\
&=&|D_1^{(N)}|+(0,0)_2^{(N)}-(p-1)(0,0)_2^{(q)}-\frac{(p-2)(q-1)}{2}\nonumber\\
&=&\frac{q+1}{4}.\nonumber
\end{eqnarray}
So we have
\begin{eqnarray}
&&\mathrm{det}(A)\nonumber\\
&=&\left\{
\begin{array}{ll}
(p+1)\left(\frac{q-1}{2}\right)^{p}\left(\left(p-1\right)^2\frac{q-1}{4}-p\right)^{\frac{q-1}{2}}\left(\frac{q-1}{4}\right)^{\frac{(p-1)(q-1)}{2}},
&\ \textrm{if $q\equiv 1\pmod{4}$;}\nonumber\\
(p+1)\left(\frac{q-1}{2}\right)^{p}\left(\left(p-1\right)^2\frac{q+1}{4}+p\right)^{\frac{q-1}{2}}\left(\frac{q+1}{4}\right)^{\frac{(p-1)(q-1)}{2}},
&\ \textrm{if $q\equiv 3\pmod{4}$.}\nonumber\\
\end{array}
\right.\label{two cases}
\end{eqnarray}
The result follows.

\begin{theorem}
Let $p$ and $q$ be twin primes satisfying $q=p+2$. Suppose $\{s_i\}_{i=0}^{N-1}$ is the Ding-Helleseth generalized cyclotomic sequence defined by Equations (\ref{eqsi-p}). Then the $2$-adic complexity $\phi_2(s)$ of $\{s_i\}_{i=0}^{N-1}$ is
\begin{eqnarray}
\phi_2(s)=N-1.\nonumber
\end{eqnarray}
\end{theorem}
\noindent{\bf Proof.} Let $r$ be a prime factor of $2^{N}-1$ and $R=Ord_r(2)$ be the multiplicative order of 2 modulo $r$. Since $2^{N}-1\equiv 0\pmod{r}$, $2^{R}\equiv 1\pmod{r}$, therefore $R\mid N$. So we can get $R\in \{p,q,pq\}.$ By Fermat's little Theorem, we know that $2^{r-1}\equiv 1\pmod{r}$. Then $R\mid r-1$, therefore $r=kR+1$, where $k$ is a positive integer. From Lemma \ref{lmgcdA-N}, we first calculate the value of $\mathrm{gcd}(det(A),2^N-1)$ for different cases.

(1) If $q\equiv 1\pmod{4}$, $q=p+2.$

Case 1.  $R=pq, r=kpq+1$.

Since $p+1<r$, $\frac{q-1}{2}<r$, and $\frac{q-1}{4}<r$, we have $\gcd(p+1, r) =1$, $\gcd(\frac{q-1}{2}, r) =1$, and $\gcd(\frac{q-1}{4}, r) =1$.
Now we consider the relationship between $\left(p-1\right)^2\frac{q-1}{4}-p$ and $r$.
Since $q=p+2$, suppose that $$\left(p-1\right)^2\frac{q-1}{4}-p=pq\frac{p-3}{4}+\frac{p+1}{4}=tr=t(kpq+1),$$ where $t$ is a positive integer. Thus $pq(p-3)+(p+1)=4ktpq+4t$.

If $4t<pq$, we have $4t=p+1$,$4tk=p-3$, then we get $0<k<1$; if $4t>pq$, we have $4tk<p-3$, we can also get $0<k<1$, these all contradict the condition that $k$ is a positive integer.

Therefore, $\gcd(\left(p-1\right)^2\frac{q-1}{4}-p, r) =1$. In summary, we have$$\gcd(\mathrm{det}(A), 2^N-1) =1.$$

Case 2.  $R=p, r=kp+1$.

Using similar method to Case 1, we have $\gcd(\frac{q-1}{2}, r) =1$, $\gcd(\frac{q-1}{4}, r)=1$. Since $r$ and $p$ are primes, $r\neq p+1$, and we have $\gcd(p+1, r) =1$.
Now we consider the relationship between $\left(p-1\right)^2\frac{q-1}{4}-p$ and $r$.
Suppose they are not coprime integers, by $q=p+2$, we can get $$\left(p-1\right)^2\frac{q-1}{4}-p=\frac{p^3-p^2-5p+1}{4}.$$
Then $r\mid p^3-p^2-5p+1.$ Since $$p^3-p^2-5p+1=(\frac{p^2}{k}-\frac{(k+1)p}{k^2}+\frac{-5k^2+k+1}{k^3})(kp+1)+(1-\frac{-5k^2+k+1}{k^3}),$$
we can get $1-\frac{-5k^2+k+1}{k^3}=0.$ Then $k^3+5k^2-k-1=0.$
However, since $k$ is a positive integer, $k^3+5k^2-k-1>0,$ contradictorily.
Therefore, $\gcd(\left(p-1\right)^2\frac{q-1}{4}-p, r) =1$. In summary, we have$$\gcd(\mathrm{det}(A), 2^N-1) =1.$$\

Case 3.  $R=q, r=kq+1$.

Using similar method to Case 2, we have $\gcd(p+1, r) =1$, $\gcd(\frac{q-1}{2}, r) =1$ and $\gcd(\frac{q-1}{4}, r)=1$.
Similarly, suppose $\left(p-1\right)^2\frac{q-1}{4}-p$ and $r$ are not coprime integers, by $q=p+2$, we can get $$\left(p-1\right)^2\frac{q-1}{4}-p=\frac{q^3-7q^2+11q-1}{4}.$$
Then $r\mid q^3-7q^2+11q-1.$ Since $$q^3-7q^2+11q-1=(\frac{q^2}{k}-\frac{(7k+1)q}{k^2}+\frac{11k^2+7k+1}{k^3})(kq+1)+(-1-\frac{11k^2+7k+1}{k^3}),$$
we can get $1+\frac{11k^2+7k+1}{k^3}=0,$ then $$k^3+11k^2+7k+1=0.$$
Since $k$ is a positive integer, that's not possible.
Therefore, $\gcd(\left(p-1\right)^2\frac{q-1}{4}-p, r) =1$. In summary, we have$$\gcd(\mathrm{det}(A), 2^N-1) =1.$$

(2) If $q\equiv 3\pmod{4}$, $q=p+2, $

Similarly to the case $q\equiv 1\pmod{4}$, $q=p+2$, we have $\gcd(\mathrm{det}(A), 2^N-1) =1.$
By Lemma 6 and Equation (8), we can get $\mathrm{gcd}\left(S(2),2^N-1\right)=1$. Thus, the $2$-adic complexity $\Phi_{2}(s)$ of the sequence is
$$\Phi_{2}(s)=\left\lfloor\mathrm{log}_2{(2^N-1)}\right\rfloor=N-1.$$
\subsection{The $2$-adic complexity of the second class of sequences}
Define the Ding-Helleseth generalized cyclotomic sequence of order 2 with respect to primes $p$ and $q$ as
\begin{equation}\label{eqsi-q}
s_i=\left\{
\begin{array}{ll}
0, \ \mathrm{if}\ i\pmod{N}\in D_0^{(N)}\cup P\cup R,\\
1, \ \mathrm{if}\ i\pmod{N}\in D_1^{(N)}\cup Q.\\
\end{array}
\right.
\end{equation}
\begin{theorem}
Let $A$ be the circulant matrix based on the Ding-Helleseth generalized cyclotomic sequence defined by Equation (\ref{eqsi-q}). Then
\begin{eqnarray}
\mathrm{det}(A)=\left\{
\begin{array}{ll}
(p-1)^q\left(\frac{q+1}{2}\right)^{p}\left(\frac{q-1}{4}\right)^{\frac{p(q-1)}{2}},
&\ \textrm{if $q\equiv 1\pmod{4}$;}\nonumber\\
(p-1)^q\left(\frac{q+1}{2}\right)^{p}\left(\frac{q+1}{4}\right)^{\frac{p(q-1)}{2}},
&\ \textrm{if $q\equiv 3\pmod{4}$.}\nonumber\\
\end{array}
\right.\label{two cases}
\end{eqnarray}
\end{theorem}
\noindent{\bf Proof.}
By Lemma \ref{lmgcdA-N},
\begin{eqnarray}
\mathrm{det}(A)&=&\prod_{a=0}^{N-1}S(\xi_N^a)\nonumber\\
&=&\prod_{a\in R}S(\xi_N^a)\prod_{a\in P}S(\xi_N^a)\prod_{a\in Q}S(\xi_N^a)\prod_{a\in D_0^{(N)}}S(\xi_N^a)\prod_{a\in D_1^{(N)}}S(\xi_N^a).\nonumber
\end{eqnarray}
By Equations (\ref{eqS(x)}) and (\ref{eqsi-q}),
\begin{eqnarray}
S(\xi_N^a)=\sum_{i\in D_1^{(N)}}\xi_N^{ai}+\sum_{i\in Q}\xi_N^{ai}.\nonumber
\end{eqnarray}
(1)If $a\in R$, obviously,
\begin{eqnarray}
S(\xi_N^a)&=&|D_1^{(N)}|+|Q|\nonumber\\
&=&\frac{(p-1)(q+1)}{2}.\nonumber
\end{eqnarray}
Thus$$\prod_{a\in R}S(\xi_N^a)=\frac{(p+1)(q-1)}{2}.$$
(2) If $a\in P$,
\begin{itemize}
\item[1)]  when $a\in D_0^{(q)}p$, by Lemmas \ref{lmPQR}, \ref{lmchi=1}, Equations (\ref{eqD0D1}) and (\ref{eqeta-q}),
\begin{eqnarray}
S(\xi_N^a)&=&(p-1)\eta_1^{(q)}+(p-1)\nonumber\\
&=&-(p-1)\eta_0^{(q)};\nonumber
\end{eqnarray}
\item[2)]  when $a\in D_1^{(q)}p$, by Lemmas \ref{lmPQR}, \ref{lmchi=1}, Equations (\ref{eqD0D1}) and (\ref{eqeta-q}),
\begin{eqnarray}
S(\xi_N^a)&=&(p-1)\eta_0^{(q)}+(p-1)\nonumber\\
&=&-(p-1)\eta_1^{(q)}.\nonumber
\end{eqnarray}
\end{itemize}
And we have
\begin{eqnarray}
\prod_{a\in P}S(\xi_N^a)&=&\left((p-1)^2\eta_0^{(q)}\eta_1^{(q)}\right)^{\frac{q-1}{2}}.\nonumber
\end{eqnarray}
If $q\equiv 1\pmod{4}$, by Lemmas \ref{lm-1}, \ref{lmcyc-num-q}, \ref{lmchi=1} and Equation (\ref{eqeta-q}),
\begin{eqnarray}
\eta_0^{(q)}\eta_1^{(q)}&=&-\frac{q-1}{4}.\nonumber
\end{eqnarray}
If $q\equiv 3\pmod{4}$, by Lemmas \ref{lm-1}, \ref{lmcyc-num-q}, \ref{lmchi=1} and Equation (\ref{eqeta-q}),
\begin{eqnarray}
\eta_0^{(q)}\eta_1^{(q)}&=&\frac{q+1}{4}.\nonumber
\end{eqnarray}
So we get
\begin{eqnarray}
\prod_{a\in P}S(\xi_N^a)=\left\{
\begin{array}{ll}
(p-1)^{q-1}\left(\frac{q-1}{4}\right)^{\frac{q-1}{2}},
&\ \textrm{if $q\equiv 1\pmod{4}$;}\nonumber\\
(p-1)^{q-1}\left(\frac{q+1}{4}\right)^{\frac{q-1}{2}},
&\ \textrm{if $q\equiv 3\pmod{4}$.}\nonumber
\end{array}
\right.\label{two cases}
\end{eqnarray}
(3) If $a\in Q$, by Lemmas \ref{lmPQR} and \ref{lmchi=1},
\begin{eqnarray}
S(\xi_N^a)=\frac{q-1}{2}\sum_{k\in Q}\xi_N^k+\sum_{k\in Q}\xi_N^k=-\frac{q+1}{2}.\nonumber
\end{eqnarray}
Then$$\prod_{a\in Q}S(\xi_N^a)=\left(\frac{q+1}{2}\right)^{p-1}.$$
(4) If $a\in D_0^{(N)}$, by Lemmas \ref{lmPQR} and \ref{lmchi=1},
\begin{eqnarray}
S(\xi_N^a)=\eta_1^{(N)}-1=-\eta_0^{(N)}.\nonumber
\end{eqnarray}
Then$$\prod_{a\in D_0^{(N)}}S(\xi_N^a)=\left(\eta_0^{(N)}\right)^{\frac{(p-1)(q-1)}{2}}.$$
(5) Similarly, if $a\in D_1^{(N)}$, by Lemmas \ref{lmPQR} and \ref{lmchi=1},
\begin{eqnarray}
S(\xi_N^a)=\eta_0^{(N)}-1=-\eta_1^{(N)}.\nonumber
\end{eqnarray}
Then$$\prod_{a\in D_1^{(N)}}S(\xi_N^a)=\left(\eta_1^{(N)}\right)^{\frac{(p-1)(q-1)}{2}}.$$
When $q\equiv 1\pmod{4}$, $-1\in D_0^{(N)}$, by Lemmas \ref{lm-1} - \ref{lmy-x}, \ref{lmchi=1} and Equation (\ref{eqeta-n}), we have
\begin{eqnarray}
\eta_0^{(N)}\eta_1^{(N)}&=&\sum_{x\in D_0^{(N)}}\sum_{y\in D_1^{(N)}}\xi_N^{y-x}\nonumber\\
&=&-\frac{q-1}{4}.\nonumber
\end{eqnarray}
When $q\equiv 3\pmod{4}$, $-1\in D_1^{(N)}$, by Lemmas \ref{lm-1} - \ref{lmy-x}, \ref{lmchi=1} and Equation (\ref{eqeta-n}), we have
\begin{eqnarray}
\eta_0^{(N)}\eta_1^{(N)}&=&\sum_{x\in D_1^{(N)}}\sum_{y\in D_1^{(N)}}\xi_N^{y-x}\nonumber\\
&=&\frac{q+1}{4}.\nonumber
\end{eqnarray}
Then we get
\begin{eqnarray}
\mathrm{det}(A)=\left\{
\begin{array}{ll}
(p-1)^q\left(\frac{q+1}{2}\right)^{p}\left(\frac{q-1}{4}\right)^{\frac{p(q-1)}{2}},
&\ \textrm{if $q\equiv 1\pmod{4}$;}\nonumber\\
(p-1)^q\left(\frac{q+1}{2}\right)^{p}\left(\frac{q+1}{4}\right)^{\frac{p(q-1)}{2}},
&\ \textrm{if $q\equiv 3\pmod{4}$.}\nonumber
\end{array}
\right.\label{two cases}
\end{eqnarray}
The result follows.

\begin{theorem}
Let $p$ and $q$ be twin primes satisfying $q=p+2$. Suppose $\{s_i\}_{i=0}^{N-1}$ is the Ding-Helleseth generalized cyclotomic sequence defined by Equation (\ref{eqsi-q}). Then the $2$-adic complexity $\phi_2(s)$ of $\{s_i\}_{i=0}^{N-1}$ is
\begin{eqnarray}
\phi_2(s)=N-1.\nonumber
\end{eqnarray}
\end{theorem}
\noindent{\bf Proof.}
Similarly to the proof of Theorem 1, since $p-1<r$, $\frac{q+1}{2}<r$, $\frac{q+1}{4}<r$, $\frac{q-1}{4}<r$, it's easy to observe that they are all coprime with $r$.
Then the $2$-adic complexity is $$\Phi_{2}(s)=\left\lfloor\mathrm{log}_2{(2^N-1)}\right\rfloor=N-1.$$

\section{Conclusion}

In this paper, using circulant matrices, we prove that  $2$-adic complexity of two classes of known generalized cyclotomic sequences are optimal.


\end{document}